\documentclass[reqno, 11pt]{article} %{amsart}
\usepackage{amsthm,amsmath,amsfonts,amssymb,mathtools,bm,framed,fixltx2e}
\usepackage{bbm}
\usepackage{mathrsfs}
\usepackage{amscd}
\usepackage{enumerate}
\usepackage[top=4cm, bottom=3cm, left=3.2cm, right=3.2cm]{geometry}

\theoremstyle{definition}

\theoremstyle{remark}

\theoremstyle{Hypothesis}

\numberwithin{equation}{section} \numberwithin{lem}{section}
\numberwithin{thm}{section} \numberwithin{prop}{section}
\numberwithin{cor}{section} \numberwithin{rem}{section}\numberwithin{hyp}{section}
\allowdisplaybreaks
\title%[Random Vortex Method]
{The optimal system, similarity reduction and group-invariant solutions to the Fractional Porous Medium equation and Fractional Dual Porous Medium equation}

\author{Ying {Yang}
\thanks{School of Mathematics, Northwest University, Xian, 710069, People's Republic of China
    yangying1@stumail.nwu.edu.cn}
\and Li-zhen Wang
\thanks{School of Mathematics, Northwest University, Xian, 710069, People's Republic of China
 Email:wanglizhen@nwu.edu.cn}
}

\begin{document}
\maketitle
\date{}

%================================================================================================================================
\textbf{Abstract}: In this paper, we developed the group analysis to the time Fractional Porous Medium Equation (FPME) and the time Fractional Dual Porous Medium Equation (FDPME). The symmetry groups and the corresponding optimal systems of these two equations are obtained. Based on the above results, similarity reductions are performed and  some explicit group invariant solutions are constructed.  

\section{Introduction}

Fractional calculus plays a significant role in various fields such as electrochemistry, viscoelasticity, rheology, biology and physics \cite{SAO,L, RP,BDST}. Though it is a generalization of the calculus with integer order, it almost emerges at the same time as classical calculus, which can be traced back to Newton and Leibniz. In recent years, fractional calculus has drawn considerable attention because it is more reasonable to use the models with the fractional-order derivative to describe the nonlinear phenomena than the integer-order models \cite{JFM}. A multitude of methods such as variational iteration method \cite{GWE}, homotopy perturbation method \cite{JHH}, Adomian decomposition method \cite{GA}, differential transform method \cite{ZOS,XLW} and the invariant subspace method \cite{VS} have been used to solve the Fractional Partial Differential Equation (FPDE) and the study of the FPDE has made substantial progress. \\
\qquad It is acknowledged that Lie symmetry analysis is an effective and systematic method for the study of the Partial Differential Equation (PDE) \cite{BL}.
This method was first proposed by Sophus Lie in nineteenth century and further generalized by Ovsianikov \cite{LVO} and others \cite{BL,PJ,BSS}. Symmetries of ordinary differential equation have been studied by Kasatkin and his collaborator \cite{RAS1,AAK}. The time fractional linear equation under scaling transformation and its solutions have been represented in \cite{EY}. Recently, some time fractional nonlinear equations have been studied, such as fractional heat conduction equation \cite{VDTM}, fractional diffusion equation \cite{RAS2}, fractional Burgers equation and Korteweg-de Vries equation \cite{RT,GXY,HL}, fractional Harry-Dym Type equation \cite{WWSH,QR} . \\
Porous Medium equation (PME) has a wide range of applications in physics, diffusion process and engineering science \cite{JDM,JC,JRK,ERA}. PME is derived from the flow of isentropic gas which is governed by the equation of state, conservation of mass and Darcy's law \cite{VJ}. Darcy's law can be derived from the resistance relationship encountered in laminar flow in porous media \cite{BMG,CHS,IC}. So we can get the one dimensional PME $\frac{\partial u}{\partial t}=(u^r)_{xx}£¬$ by eliminating the $p$, $v$ in Darcy's law. In recent years, this type equations have aroused widespread concern. Estevez and Qu constructed functional seperation of variables solutions of the generalized PME by the generalized conditional symmetry approach \cite{PCS}. Bonforte and Grillo studied the asymptotics of the PME using Sobolev inequalities \cite{BMG}. Caffarelli and Vazquez introduced the one dimensional fractional version of the porous medium equation and proved the existence of weak and bounded solutions propagating with finite speed \cite{LCJLV}. Carrillo and Huang deduced the exponential convergence towards stationary states for the one dimensional PME with fractional pressure \cite{CHS}. Stan and Teso investigated the type of solutions by the transformations of Self-Similar solutions \cite{SDF}. Pablo and Quiros developed a theory of existence, uniqueness and regularity for the PME with fractional diffusion and constructed an $L^1$-contraction semigroup \cite{AFAJ}. It were Biler and Imbert who showed the existence of sign changing weak solutions and constructed explicit compactly supported self-solutions of nonlocal PME \cite{BPCKG}. And the Dual Porous medium equation has been discussed in \cite{VS} with the invariant subspace method.\\
For the sake of obtaining more information to comprehend the Porous Media type equations, by the Lie symmetry analysis, we consider the following one dimensional time FPME and time Fractional Dual Porous Medium equation (FDPME) which read:
\begin{equation}\label{PME}
\partial_t^\alpha u =(u^r)_{xx}=r(r-1)u^{r-2}u_x^2+ru^{r-1}u_{xx},  \\
\end{equation}
and
\begin{equation}\label{DPME0}
\partial_t^\alpha u=a(u_{xx})^2+bu_xu_{xx}+c[uu_{xx}-\frac{2}{3}(u_x)^2],\\
\end{equation}
respectively, where $u=u(x,t)$ is the unknown and represents the scaled density. $0<\alpha<1, r>0$, $r\neq1$ and $a,b,c$ are real numbers, $\partial_t^\alpha$ is the Riemann-Liouville fractional derivative operator with order $\alpha$.  \\
This paper is organized as follows. In section 2, we introduce the definition of Riemann-Liouville derivative with related properties and the approach of Lie symmetry analysis. In Section 3, we use the Lie symmetry analysis to solve the time FPME. In addition, we construct the optimal system of subgroups admitted by the FPME. Also, the similarity reduction is performed and group-invariant solutions of FPME are obtained. In section 4, the symmetry algebra admitted by FDPME and the corresponding optimal system are derived. Furthermore, the group-invariant solutions to FDPME is established after finishing the similarity reduction  .

%================================================================================================================================
\section{Notation and Preliminaries}
\textbf{Definition2.1}: Suppose $n\in N$. The Riemann-Liouville fractional derivative is defined by
\begin {equation}\label{RL derivative}
\partial_t^\alpha u=
\left\{
  \begin{aligned}
&\frac{\partial^n u}{\partial t^n},          ~~\mbox{if}~~\alpha=n,\\
&\frac{1}{\Gamma(n-\alpha)}\frac{\partial^n}{\partial t^n}\int_0^t(t-s)^{n-\alpha-1}u(t,s)ds,~~\mbox{if}~~0<n-1<\alpha<n.\\
  \end{aligned}
    \right.
\end{equation}
Where $\frac{\partial^n}{\partial t^n}$ is the usual partial derivative of integer order $n$ and $\Gamma(x)=\int_0^\infty t^{x-1}e^{-t}dt$.
Some useful formulas and properties of Riemann-Liouville derivative were summarized in \cite{JFM} and one of them is listed as follows:
\begin{equation}\label{mRL properties1}
\partial^\alpha_t t^\gamma=\frac{\Gamma(\gamma+1)}{\Gamma(\gamma+1-\alpha)} t^{\gamma-\alpha},\gamma>0,\\
\end{equation}
Consider the general time FPDE with the form
\begin{equation}\label{FPDE}
\partial_t^\alpha u=F(t,x,u_x,u_{xx},\cdots).\\
\end{equation}
Assume \eqref{FPDE} is invariant under a one parameter Lie group of point transformations
\begin{equation}\label{point transformation}
 \begin{aligned}
 t^*=t+\epsilon\tau(x,t,u)+O(\epsilon^2),\\
 x^*=x+\epsilon\xi(x,t,u)+O(\epsilon^2),\\
 u^*=u+\epsilon\eta(x,t,u)+O(\epsilon^2),\\
\partial^\alpha_{t^*}{u^*}=\partial_t^\alpha u+\epsilon\eta_\alpha^0(x,t,u)+O(\epsilon^2),\\
\partial_{x^*}{u^*}=\partial_x u+\epsilon\eta^x(x,t,u)+O(\epsilon^2),\\
\partial_{x^*}^2{u^*}=\partial_x^2 u+\epsilon\eta^{xx}(x,t,u)+O(\epsilon^2).\\
 \end{aligned}
\end{equation}
Where $\tau$, $\xi$, $\eta$ are infinitesimals and $\eta_\alpha^0$, $\eta^x$, $\eta^{xx}$ are extended infinitesimal of orders $\alpha$, $1$ and $2$, respectively. The formulas of $\eta_\alpha^0$, $\eta^x$, $\eta^{xx}$ are
\begin{equation}\label{pt11}
 \begin{aligned}
\eta_\alpha^0&=D_t^\alpha(\eta)+\xi D_t^\alpha(u_x)-D_t^\alpha(\xi u_x)+D_t^\alpha(uD_t(\tau))-D_t^{\alpha+1}(\tau u)+\tau D_t^{\alpha+1}(u)\\
&=\frac{\partial^\alpha \eta}{\partial t^\alpha}+(\eta_u-\alpha D_t{\tau})\frac{\partial^\alpha u}{\partial t^\alpha}-u\frac{\partial^\alpha \eta_u}{\partial t^\alpha}+\mu+\sum_{n=1}^\infty[\mathrm{C}_\alpha^n\frac{\partial^n \eta_u}{\partial t^n}
-\mathrm{C}_\alpha^{n+1}D_t^{n+1}(\tau)]D_t^{\alpha-n}(u)\\
&-\sum_{n=1}^{\infty}\mathrm{C}_\alpha^nD_t^n(\xi)D_t^n(\xi)D_t^{\alpha-n}(u_x).
  \end{aligned}
\end{equation}
where$$\mathrm{C}_\alpha^n=\frac{\Gamma(\alpha+1)}{\Gamma(n+1)\Gamma(\alpha+1-n)},$$
$$\mu=\sum_{n=2}^\infty\sum_{m=2}^n\sum_{k=2}^m\sum_{r=0}^{k-1}\mathrm{C}_\alpha^n\mathrm{C}_n^m\mathrm{C}_k^r\frac{1}{k!}\frac{t^{n-\alpha}}{\Gamma(n+1-\alpha}[-u]^r\frac{\partial^m}{\partial t^m}[u^{k-r}]\frac{\partial^{n-m+k}\eta}{\partial t^{n-m}\partial u^k},$$
\begin{equation}\label{pt12}
 \eta^x=D_x(\eta)-u_tD_x(\tau)-u_xD_x(\xi),\\
\end{equation}
\begin{equation}\label{pt13}
 \eta^{xx}=D_x(\eta^x)-u_{xt}D_x(\tau)-u_{xx}D_x(\xi),\\
\end{equation}
and the symbols $D_t$ and $D_x$ indicate the total derivatives with respect to $t$ and $x$ respectively, that is
\begin{equation*}
D_t=\partial_t+u_t\partial_u+u_{tt}\partial_{u_t}+u_{xt}\partial{u_x}+\cdots,
\end{equation*}
\begin{equation*}
D_x=\partial_x+u_x\partial_u+u_{tx}\partial_{u_t}+u_{xx}\partial{u_x}+\cdots,
\end{equation*}
with the infinitesimal generator
\begin{equation}\label{infinitesimal generator}
V=\tau\partial_t+\xi\partial_x+\eta\partial_u,\\
\end{equation}
where
$$\tau=\frac{dt^*}{d\epsilon}|_{\epsilon=0}, \mbox{~~~~} \xi=\frac{dx^*}{d\epsilon}|_{\epsilon=0}, \mbox{~~~~} \eta=\frac{du^*}{d\epsilon}|_{\epsilon=0}.$$
According to the infinitesimal invariance criterion, \eqref{FPDE} admits \eqref{point transformation} if and only if the prolonged vector field  $pr^{(\alpha,2)}V$ annihilates \eqref{FPDE} on its solution manifold \cite{PJ}, namely,
\begin{equation*}\label{vector field}
pr^{(\alpha,2)}V(\triangle_1)|_{\triangle_1=\partial_t^\alpha u-F=0}=0,\\
 \end{equation*}
And the operator $pr^{(\alpha,2)}V$ takes the form
\begin{equation}
pr^{(\alpha,2)}V=V+\eta_\alpha^0\partial_{\partial_t^\alpha u}+\eta^x\partial_{u_x}+\eta^{xx}\partial_{u_{xx}},\\
\end{equation}
where $\eta_\alpha^0, \eta^x, \eta^{xx}$ satisfy \eqref{pt11}, \eqref{pt12}, \eqref{pt13}, respectively.\\
Since we intend to construct the group-invariant solution, the definition of invariant solutions is given as follows.\\
\textbf{Definition 2.2} $u=\Theta(x,t)$ is an invariant solution of \eqref{FPDE} corresponding to \eqref{infinitesimal generator} if and only if
(i) $ u=\Theta (x,t)$ is an invariant surface of \eqref{infinitesimal generator}, i.e
$$\xi(x,t,\Theta)\Theta_x+\tau(x,t,\Theta)\Theta_t=\eta(x,t,\Theta).$$
(ii) $ u=\Theta (x,t)$ solves \eqref{FPDE}.\\
We can find the following results from \cite{AHS}, which will help us solve equations later.\\
\textbf{Lemma 2.1}\cite{AHS} (i) Let $\alpha>0$. If $\frac{\beta+\alpha}{1-r}>-1$, $\lambda,\beta\in R$, then the $\alpha$-th order fractional differential equation
\begin{equation}\label{case 1}
\partial^\alpha_t g(t)=\lambda (t-a)^\beta[g(t)]^r, (t>a; m>0 , m\neq 1),\\
\end{equation}
 has the following explicit solutions:
\begin{equation}\label{case 1.1}
g(t)=[\frac{\Gamma(\frac{\alpha+\beta}{r-1}+1)}{\lambda\Gamma(\frac{\alpha r+\beta}{r-1}+1)}]^{\frac{1}{r-1}} t^{\frac{\alpha+\beta}{1-r}}.\\
\end{equation}
(ii) Let $\alpha>0$. If $2(\alpha+\beta)>-1$, then the $\alpha$-th order fractional differential equation
\begin{equation}\label{case 2}
\partial^\alpha_t g(t)=\lambda (t-a)^\beta[g(t)]^\frac{1}{2},  (t>a; \lambda, \beta\in R, \lambda\neq0),    \\
\end{equation}
has the exact solution
\begin{equation}\label{case 21}
g(t)=[\frac{\lambda\Gamma(\alpha+2\beta+1)}{\Gamma(2\alpha+2\beta+1)}]^2(t-a)^{2(\alpha+\beta)}.\\
\end{equation}
(iii) Let $\alpha>0$. If $\alpha+\beta<1$, then the $\alpha$-th order fractional differential equation:
\begin{equation}\label{case 3}
\partial^\alpha_t g(t)=\lambda (t-a)^\beta[g(t)]^2,   (t>a;\lambda, \beta\in R, \lambda\neq0), \\
\end{equation}
has the exact solution
\begin{equation}\label{case 31}
g(t)=\frac{\Gamma(1-\alpha-\beta)}{\lambda\Gamma(1-2\alpha-\beta)}(t-a)^{-(\alpha+\beta)}.\\
\end{equation}

\section{Lie symmetry analysis, optimal system and similarity reduction for the time FPME }
In this section, we will employ the symmetry analysis method introduced above to construct the symmetry group of the time FPME. Then we establish the optimal system of the obtained algebra and the group invariant solutions.
\subsection{Lie symmetry analysis for the time FPME}
For the time FPME \eqref{PME}, the invariance criterion takes the form
\begin{equation*}
[\eta_\alpha^0-2r(r-1)u^{r-2}u_x\eta^x-ru^{r-1}\eta^{xx}-r(r-1)(r-2)u^{r-3}\eta u_x^2-r(r-1)u^{r-2}\eta u_{xx}]|_{\partial_t^\alpha u ={(u^r)}_{xx}}=0.
\end{equation*}
Taking the coefficients of linearly independent derivatives $u_t$, $u_x$, $u_{xx}$, $\cdots$, and $\partial_t^{\alpha-n} u_x$, $\partial_t^{\alpha-n} u$ to be zero in the above relation, we can obtain the following determining equations.
\begin{equation*}
 \begin{aligned}
             &\mathrm{C}_\alpha^n\partial_t^n (\eta_u)-\mathrm{C}_\alpha^{n+1}D_t^{n+1}(\tau)=0 ,      n=1,2,3,...,\\
             &(r-1)\eta_x+u(\eta_{xu}-\xi_{xx})=0,\\
             &\alpha D_t{(\tau)}u-2u\xi_x+(r-1)\eta=0,\\
  &\alpha(r-1)u\tau_t+(r-1)u\eta_u-2(r-1)u\xi_x+u^2\eta{uu}-2u^2\xi_{xu}+(r-1)(r-2)\eta=0,\\
 &\tau_x=\tau_u=\xi_t=\xi_u=\eta_{xx}=0,
 \end{aligned}
\end{equation*}
Solving the above determining equations yields the following general solutions:
\begin{equation}\label{111}
\tau=C_1t+C_2,\mbox{~~~~}   \xi=C_3x+C_4,\mbox{~~~~}   \eta=\frac{2C_3-\alpha C_1}{r-1}u.\\
\end{equation}
At the same time, in order to satisfy the invariant condition of Lie symmetry for fractional differential equation \eqref{111}, we have
$$\tau(x,t,u)|_{t=0}=0,$$
Thus, the infinitesimal generator of \eqref{PME} is
\begin{equation*}
V=C_1t\partial_t+(C_3x+C_4)\partial_x+(\frac{2C_3-\alpha C_1}{r-1})u\partial_u.\\
\end{equation*}
Therefore, we can obtain the following result on the symmetry group of equation \eqref{PME}:\\
\textbf{Theorem 3.1} The symmetry group of the time FPME \eqref{PME} is spanned by the following vector fields
\begin{equation}\label{VECTOR}
V_{11}=t\partial_t-\frac{\alpha}{r-1}u\partial_u,\mbox{~~~~}
V_{12}=x\partial_x+\frac{2}{r-1}u\partial_u,\mbox{~~~~}
V_{13}=\partial_x.\\
\end{equation}
\subsection{Optimal system for algebra \eqref{VECTOR} of the time FPME }
In order to find group-invariant solutions of \eqref{PME}, we need establish the optimal system of \eqref{VECTOR} first. The nonzero commutators of the algebra \eqref{VECTOR} read:
\begin{equation}\label{CR}
 [V_{12},V_{13}]=-V_{13},\mbox{~~~~} [V_{13},V_{12}]=V_{13}.\\
\end{equation}
The action of the adjoint operator is given by the following Lie series
\begin{equation}\label{adjoint}
\mathrm{Ad}(\exp(\epsilon V_i))V_j=V_j-\epsilon [V_i,V_j]+\frac{\epsilon^2}{2}[V_i,[V_i,V_j]]-\cdots,\\
\end{equation}
where $[V_{i},V_{j}]$ $(i,j=1,2,3)$ is the commutator and $\epsilon$ is a parameter.
By the commutation relationship \eqref{CR} and the formula \eqref{adjoint}, we can get the adjoint action listed in Table 1.

\renewcommand{\thetable}{\arabic{table}}
\renewcommand{\arraystretch}{1}
\begin{table}[h]
\caption{{\small  The adjoint representation of $H_1$ on $h_1$.}}
$$\begin{array}{ccccccc}\hline
\mbox{$Ad(\varepsilon\cdot)$ }&\mbox{$V_{11}$}&\mbox{$V_{12}$}&\mbox{$V_{13}$}\\[.04cm]\hline
V_{11}&V_{11}&V_{12}&V_{13}\\
V_{12}&V_{11}&V_{12}&e^\varepsilon V_{13}\\
V_{13}&V_{11}&V_{12}-\varepsilon V_{13}&V_{13}\\\hline
\end{array}$$
\end{table}
Following the approach of Ovsiannikov \cite{LVO}, we obtain the one-dimensional optimal system of the algebra \eqref{VECTOR}.\\
\textbf{Theorem 3.2} The one-dimensional optimal system of \eqref{VECTOR} is given by
\begin{equation}\label{OS1}
\begin{aligned}
&r_{11}=V_{11},\mbox{~} r_{12}=V_{12}, \mbox{~}r_{13}=V_{13},\mbox{~}  r_{14}=V_{11}+V_{13}, \mbox{~}r_{15}=V_{11}-V_{13}, \mbox{~}r_{16}=V_{11}+\gamma V_{12}, (\gamma\neq0).
\end{aligned}
\end{equation}
\subsection{Similarity reduction and the group invariant solutions of time FPME}
In this subsection, according to the optimal system obtained above, we intend to construct the group invariant solutions to FPME \eqref{PME}.
For each inequivalent subalgebra in Theorem 3.1, we perform the similarity reduction and deduce the reduced nonlinear fractional determining equations and construct the corresponding group-invariant solutions to \eqref{PME}. \\
Case 1: $r_{11}=t\partial_t-\frac{\alpha}{r-1}u\partial_u.$ Integrating the invariant surface condition $ \frac{dt}{t} =\frac{du}{-\frac{\alpha}{r-1}u }=\frac{dx}{0},$ we can get the invariant solution \eqref{PME} with the form
\begin{equation}\label{invariant solution 1}
u(t,x)=t^{-\frac{\alpha}{r-1}}f(x).\\
\end{equation}
In order to find the explicit expression of the function $f(x)$, substituting \eqref{invariant solution 1} into \eqref{PME}, we have
\begin{equation}\label{y1}
f(x)\partial_t^\alpha t^{-\frac{\alpha}{r-1}}=rt^{-\frac{\alpha r}{r-1}}f^{r-2}(x)[(r-1)f'^2(x)+f(x)f''(x)].\\
\end{equation}
 Applying the formula \eqref{mRL properties1} to \eqref{y1} yields
\begin{equation}\label{DE1}
r\frac{\Gamma(\alpha-\frac{\alpha r}{r-1}+1)}{\Gamma(\frac{\alpha r+1}{1-r})}f^{r-2}(x)[(r-1)f'(x)^2+f(x)f''(x)]=f(x).\\
\end{equation}
Therefore, the solution of time FPME \eqref{PME} has the following form
\begin{equation*}
u(t,x)=t^{-\frac{\alpha}{r-1}}f(x),\\
\end{equation*}
where $f(x)$ is determined by \eqref{DE1}.\\
Case 2: $r_{12}=x\partial_x+\frac{2}{r-1}u\partial_u.$ According to the vector field $V_2$, we can get the invariant solution of time FPME \eqref{PME} with the form
\begin{equation}\label{invariant solution 2}
u(t,x)=x^{\frac{2}{r-1}}g(t).\\
\end{equation}
Substituting \eqref{invariant solution 2} into \eqref{PME}, we obtained
\begin{equation}\label{invariant solution 212}
\partial_t^\alpha g=\frac{2r(r+1)}{(r-1)^2}g^r(t).\\
\end{equation}
Thus, with the help of Lemma 2.1, we can find some special solutions to \eqref{PME}.\\
\textbf{Theorem 3.3} The time FPME \eqref{PME} has the following solutions:\\
(i) When $\frac{\alpha}{1-r}>-1,$ solutions of \eqref{PME} are:
\begin{equation*}
u(t,x)=x^{\frac{2}{r-1}}[\frac{(r-1)^2}{2r(r+1)}\frac{\Gamma(\frac{\alpha}{r-1}+1)}{\Gamma(\frac{\alpha r}{r-1}+1))}]^{\frac{1}{r-1}} t^{\frac{\alpha}{1-r}}.\\
\end{equation*}
(ii) When $r=\frac{1}{2}$, $\alpha>-\frac{1}{2}$,
solution of \eqref{PME} is:$$u(t,x)=x^{-4}[6\frac{\Gamma(\alpha+1)}{\Gamma(2\alpha+1)}]^2 t^{2\alpha}.$$
(iii) When $r=2$, $\alpha<\frac{1}{2}$,
solution of \eqref{PME} is:$$u(t,x)=\frac{1}{12}x^2\frac{\Gamma(1-\alpha)}{\Gamma(1-2\alpha)} t^{-\alpha}.$$
\textbf{Proof:} (i) Take $\lambda=\frac{2r(r+1)}{(r-1)^2}$, $\beta=0$ into \eqref{case 1.1}. In Lemma 2.1 (i) for $\frac{\alpha}{1-r}>-1,$ we can get the solution of equation \eqref{PME}
\begin{equation*}
g(t)=[\frac{(r-1)^2}{2r(r+1)}\frac{\Gamma(\frac{\alpha}{r-1}+1)}{\Gamma(\frac{\alpha r}{r-1}+1))}]^{\frac{1}{r-1}} t^{\frac{\alpha}{1-r}}.\\
\end{equation*}
Putting the above equation into \eqref{invariant solution 2}, we verify case 1.\\
(ii) Substituting $a=0$, $\beta=0$, $r=\frac{1}{2}$, $\lambda=6$ into \eqref{case 2} in Lemma 2.1 (ii). For $\alpha>-\frac{1}{2}$, equation \eqref{invariant solution 212} has exact solution
\begin{equation*}
g(t)=[6\frac{\Gamma(\alpha+1)}{\Gamma(2\alpha+1)}]^2 t^{2\alpha}.
\end{equation*}
And the solution of equation \eqref{invariant solution 2} has the exact form:$$u(t,x)=x^{-4}[6\frac{\Gamma(\alpha+1)}{\Gamma(2\alpha+1)}]^2 t^{2\alpha}.$$
(iii) In Lemma 2.1 (iii), substituing $a=0$, $\beta=0$, $r=2$, $\lambda=12$ into \eqref{case 3}, we can get that for $\alpha<1$, the equation \eqref{invariant solution 212} has the exact solution
\begin{equation*}
g(t)=\frac{1}{12}\frac{\Gamma(1-\alpha)}{\Gamma(1-2\alpha)} t^{-\alpha}.\\
\end{equation*}
Then the equation \eqref{invariant solution 2} has the exact form:$$u(t,x)=\frac{1}{12}x^2[\frac{\Gamma(1-\alpha)}{\Gamma(1-2\alpha)}]^2 t^{-\alpha}.$$
Case 3: $r_{13}=V_{13}=\partial_x$. According to the invariant surface condition, we get the similarity variables $t$, $u$. Thus we have $u=g(t)$. Putting it into \eqref{invariant solution 212} yields the fractional ODE
$$\partial_t^\alpha g(t)=0,$$
which implies that for $\lambda\in R$, $g(t)=\lambda t^{\alpha-1}$ is a solution of \eqref{invariant solution 212}.
  So, the equation \eqref{PME} has the solutions of the following form
\begin{equation}
u(t,x)=\lambda t^{\alpha-1}.
\end{equation}
Case 4: $r_{14}=V_{11}+V_{13}=t\partial_t-\frac{\alpha u}{r-1}\partial_u+\partial_x.$
Since the characteristic equation is $$\frac{dt}{t} =\frac{du}{\frac{-\alpha u}{r-1}}=\frac{dx}{1},$$
we can get the similarity variables $e^{\frac{\alpha}{r-1}x}u$ and $te^{-x}$. Denote $\beta=r\frac{\alpha r+r-1}{r-1}.$ The group invariant solution has the form
$$u(t,x)=f(te^{-x})e^{-\frac{\alpha}{r-1}x},$$
where $f(z)(z=te^{-x})$ satisfies:
$$\partial_z^\alpha f(z)=\frac{(\alpha r)^2}{(r-1)^2}f^r+\frac{\alpha r^2}{r-1}zf^{r-1}f'+r(r-1)z^2f^{r-2}f'^2+rz^2f^{r-1}f''(z)+\beta zf^{r-1}f'.$$
Case 5: $r_{15}=V_{11}-V_{13}=t\partial_t-\frac{\alpha}{r-1}u\partial_u-\partial_x.$ We arrived at the two similarity variables $te^x$ and $ue^{-\frac{\alpha}{r-1}x}$. The group-invariant solution has the form
$$u=f(te^x)e^{-\frac{\alpha}{r-1}x},$$
where $f(z)(z=te^x)$ is determined by
\begin{equation*}
\partial_z^\alpha f(z)=\frac{(\alpha r)^2}{(r-1)^2}f^r+\frac{\alpha r^2}{r-1}zf^{r-1}f'+r(r-1)z^2f^{r-2}f'^2+rz^2f^{r-1}f''+\beta zf^{r-1}f'.
\end{equation*}
Case 6: $r_{16}=V_{11}+\gamma V_{12}=t\partial_t+\gamma x\partial_x+\frac{2\gamma-\alpha}{r-1}u\partial_u.$ Since the similarity variables are $tx^{-\frac{1}{\gamma}}$ and $ux^{-\frac{2\gamma-\alpha}{\gamma(r-1)}}$, we get the solutions with the following form
$$u=f(tx^{-\frac{1}{\gamma}})x^{\frac{2\gamma-\alpha}{\gamma(r-1)}},$$
where $f(z)(z=tx^{-\frac{1}{\gamma}})$ satisfies the following fractional ODE
\begin{equation*}
 \begin{aligned}
  &\partial_z^\alpha f(z)=\frac{(2\gamma-\alpha)r}{\gamma(r-1)}\frac{\gamma r-\alpha r+\gamma}{\gamma(r-1)}f^r(z)-\frac{(2\gamma-\alpha)r^2}{\gamma^2(r-1)}zf^{r-1}(z)f'(z)+\frac{r(r-1)}{\gamma^2}z^2f^{r-2}(z)f'^2(z)\\
  &+\frac{r}{\gamma^2}z^2f^{r-1}(z)f''(z)
  -\frac{r(\gamma r-\alpha r-r+1+\gamma)}{\gamma(r-1)}zf^{r-1}(z)f'(z).\\
 \end{aligned}
 \end{equation*}
\section{Lie symmetry analysis, optimal system and group-invariant solutions for the time FDPME}
In this section, we study the dual porous medium equation by the Lie symmetry analysis.
\subsection{Lie symmetry for the time FDPME}
For the time FDPME \eqref{DPME0}, the invariance criterion takes the form
\begin{equation*}
\eta_\alpha^0-2au_{xx}\eta^{xx}-b\eta^xu_{xx}-bu_x\eta^{xx}-c\eta u_{xx}-cu\eta^{xx}+\frac{4}{3}cu_x\eta^x=0.
\end{equation*}
The determining equations of \eqref{DPME0} are
\begin{equation*}
 \begin{aligned}
 &-b\eta_{xx}-cu(2\eta_{xu}-\xi_{xx})+\frac{4}{3}c\eta_x=0,\\
 &\frac{2}{3}c\eta_u+\frac{2}{3}c\alpha D_t(\tau)-2b\eta_{xu}+b\xi_{xx}-cu\eta_{uu}+2cu\xi_{xu}-\frac{4}{3}c\xi_x=0,\\
 &-\alpha cuD_t(\tau)-2a\eta_{xx}-b\eta_x-c\eta+2cu\xi_x=0,\\
\end{aligned}
\end{equation*}
\begin{equation*}
 \begin{aligned}
 &-a\eta_u+4a\xi_x-aD_t(\tau)=0,\\
 &-\alpha bD_t(\tau)-4a\eta_{xu}+2a\xi_{xx}-b\eta_u+3b\xi_x=0,\\
 &\mathrm{C}_\alpha^n\partial_t^n (\eta_u)-\mathrm{C}_\alpha^{n+1}D_t^{n+1}(\tau)=0 ,      n=1,2,3,...,\\
 &\tau_x=\tau_u=\xi_t=\xi_u=\eta_{uu}=0,
 \end{aligned}
\end{equation*}
Solving the above determining equations, we obtained
$$\tau=(D_1t+D_2)\partial_t,\mbox{~~~~}\xi=D_3\partial_x,\mbox{~~~~}\eta=(-\alpha D_1u+D_4)\partial_u.$$
According to the invariant condition of Lie symmetry, we find that the infinitesimal generator of \eqref{DPME0} is $$W=D_1t\partial_t+D_3\partial_x+(-\alpha D_1u+D_4)\partial_u.$$
Therefore, we can obtain the following result on the symmetry group of equation \eqref{DPME0}.\\
\textbf{Theorem 4.1} The symmetry group of the time FDPME \eqref{DPME0} is spanned by the following vector fields
\begin{equation}\label{VECTOR2}
V_{21}=t\partial_t-\alpha u\partial_u,\mbox{~~~~}
V_{22}=\partial_x,\mbox{~~~~}
V_{23}=\partial_u.\\
\end{equation}
\subsection{Optimal system and the similarity reduction for the FDPME}
Similarly, in order to find group-invariant solutions of \eqref{DPME0}, we need establish the optimal system of \eqref{VECTOR2} at first. The nonzero commutators of the algebra \eqref{VECTOR2} are
\begin{equation}\label{CR2}
 [V_{21},V_{23}]=\alpha V_{23},\mbox{~~~~} [V_{23},V_{21}]=-\alpha V_{23}.\\
\end{equation}
By the commutation relationship \eqref{CR2} and the formula \eqref{adjoint}, we can get the adjoint action listed in Table 2.

\renewcommand{\thetable}{\arabic{table}}
\renewcommand{\arraystretch}{1}
\begin{table}[h]
\caption{{\small  The adjoint representation of $H_2$ on $h_2$.}}
$$\begin{array}{ccccccc}\hline
\mbox{$Ad(\varepsilon\cdot)$ }&\mbox{$V_{21}$}&\mbox{$V_{22}$}&\mbox{$V_{23}$}\\[.04cm]\hline
V_{21}&V_{21}&V_{22}&e^{-\alpha\varepsilon} V_{23}\\
V_{22}&V_{21}&V_{22}& V_{23}\\
V_{23}&V_{21}-\alpha\varepsilon V_{23}&V_{22}&V_{23}\\\hline
\end{array}$$
\end{table}
And the optimal system of \eqref{VECTOR2} follows from the method introduced in \cite{JFM}.\\
\textbf{Theorem4.2} The one-dimensional optimal system of \eqref{VECTOR2} is given by
\begin{equation}\label{OS2}
\begin{aligned}
&r_{21}=V_{21}, r_{22}=V_{22}, r_{23}=V_{23}, r_{24}=V_{22}+V_{23}, r_{25}=V_{22}-V_{23}, r_{26}=V_{21}+\rho V_{22},(\rho\neq0).
\end{aligned}
\end{equation}
With the help of the obtained algebra, we can construct the group-invariant solutions to equation \eqref{DPME0} as follows.\\
case 1: $r_{21}=V_{21}=t\partial_t-\alpha u\partial_u.$
The characteristic equation of $r_{21}$ is $$\frac{dt}{t}=\frac{du}{-\alpha u}=\frac{dx}{0}.$$
And the invariance of $r_{21}$ are $ut^\alpha$ and $x$. Thus the group-invariant solution of \eqref{DPME0} is $$u=t^{-\alpha}f(x),$$
where for $0<\alpha<\frac{1}{2}$, $f(x)$ is defined by
$$\frac{\Gamma(1-\alpha)}{(1-2\alpha)}f(x)=af''^2(x)+bf'(x)f''(x)+cf''(x)-\frac{2}{3}cf'^2(x). $$
case 2: $r_{22}=V_{22}=\partial_x.$
Since the characteristic equation of $r_{22}$ is $$\frac{dt}{0}=\frac{du}{0}=\frac{dx}{1},$$
we obtained the group-invariant solution of \eqref{DPME0} with the form: $u=g(t),$
where $g(t)$ satisfies $\partial_t^\alpha g(t)=0.$
Then we get $g(t)=\kappa t^{\alpha-1}$. So, for $\kappa \in R$, $$u=\kappa t^{\alpha-1}.$$
case 3: $r_{24}=V_{22}+V_{23}=\partial_x+\partial_u.$
Since the characteristic equation of $r_{24}$ is $$\frac{dx}{1}=\frac{dt}{0}=\frac{du}{1},$$
the invariant variables is $u-x$ and $t$ and the invariant solution of \eqref{DPME0} is $$u(t,x)=x+f(t)$$
where $\partial_t^\alpha f(t)=-\frac{2}{3}c.$\\
case 4: $r_{25}=V_{22}-V_{23}=\partial_x-\partial_u.$
The characteristic equation of $r_{25}$ is $$\frac{dx}{1}=\frac{dt}{0}=\frac{du}{-1}.$$
The invariant variables are $u+x$ and $t$. And the invariant solution of \eqref{DPME0} is $$u(x,t)=-x+f(t),$$
where $\partial_t^\alpha f(t)=-\frac{2}{3}c.$\\
case 5: $r_{26}=V_{21}+\rho V_{22}=t\partial_t-\alpha u\partial_u+\rho\partial_x.$
According to the characteristic equation of $r_{26}$ $$\frac{dt}{t}=\frac{du}{-\alpha u}=\frac{dx}{\rho},$$
we arrive at the invariant variables $te^{-\frac{x}{\rho}}$ and $ue^{-\frac{\alpha}{\rho}x}$, and the invariant solution of \eqref{DPME0}  $$u(t,x)=f(te^{-\frac{x}{\rho} })e^{-\frac{\alpha}{\rho}x},$$
where $f(z)(z=te^{-\frac{x}{\rho}})$ satisfies
\begin{equation}
\begin{aligned}
\partial_z^\alpha f(z)&=a(\frac{1+2\alpha}{\rho}zf'(z)+\frac{1}{\rho^2}f''(z)+\frac{\alpha^2}{\rho^2}f(z))^2-\frac{2}{3}c(zf'(z)+f(z))^2\\
&+c(\frac{1+\alpha}{\rho^2}zf''(z)f(z)+\frac{1}{\rho^2}z^2f(z)f''(z)+\frac{\alpha^2}{\rho^2}f^2(z)+\frac{\alpha}{\rho^2}zf'(z)f(z))\\
&-b(\frac{1+\alpha}{\rho^2}z^2f'^2(z)+\frac{1}{\rho^3}z^3f'(z)f''(z)+\frac{\alpha^2}{\rho^3}zf(z)f'(z)
+\frac{\alpha}{\rho^3}z^2f'^2(z)\\
&+\frac{\alpha(1+\alpha)}{\rho^3}zf'(z)f(z)+\frac{\alpha}{\rho^3}z^2f''(z)f(z)+\frac{\alpha^3}{\rho^3}f^2(z)+\frac{\alpha^2}{\rho^3}zf(z)f'(z)).\\
\end{aligned}
\end{equation}
\section{Concluding remarks}
In this paper, we investigated the Lie algebra, the optimal system, similarity reductions, and group-invariant solutions for the time FPME and time FDPME.
Through the similarity reductions, some group-invariant solutions which are determined by fractional ODEs have been obtained.
We perform elaborate analysis of subalgebra structure admitted by time FPME and time DPME with Riemann-Liouville derivative. Through investgating, we find that Lie symmetry is an efficient way for solving the FPDEs.
Solving the reduced fractional ODEs is a challenging problem. In order to obtain more explicit solutions, one need to develop some approach to overcome it. This work will be done in the future.

\textbf{Acknowledgments}

Supported by the National Natural Science Foundation of China (Grant No. 11771352, 11871396), the National Natural Science Foundation of Shaanxi Province(Grant No. 2018JM1005).

%=============================================================================================================================================
%\bibliography{fractional}
%\bibliographystyle{abbrv}

\end{document}